\documentclass[a4paper,12pt]{article}
\usepackage{amsthm,amsmath}
\usepackage{amssymb}
\usepackage{amsfonts}
\usepackage[latin1]{inputenc}
\usepackage[all]{xy}
\usepackage{enumerate}
\usepackage{a4wide}
\usepackage{latexsym}
\usepackage{mathrsfs}

\swapnumbers
\newtheorem{theorem}{Theorem}[section]
\newtheorem{proposition}[theorem]{Proposition}
\newtheorem{corollary}[theorem]{Corollary}
\newtheorem{lemma}[theorem]{Lemma}
\newtheorem*{theorem*}{Theorem}
\newtheorem*{proposition*}{Proposition}
\newtheorem*{corollary*}{Corollary}
\newtheorem*{lemma*}{Lemma}
\theoremstyle{definition}

\newtheorem{example}[theorem]{Example}
\newtheorem{remark}[theorem]{Remark}
\newtheorem*{remark*}{Remark}

\newtheorem*{definition*}{Definition}
\numberwithin{equation}{section}

\newcommand{\cat}[1]{\mathcal{#1}}
\newcommand{\coring}[1]{\mathfrak{#1}}
\newcommand{\tensor}[1]{\otimes_{#1}}

\newcommand{\rcomod}[1]{\mathsf{Comod}_{#1}}
\newcommand{\rmod}[1]{\mathsf{Mod}_{#1}}

\renewcommand{\hom}[3]{\mathrm{Hom}_{#1}\left( #2 \, , \, #3\right)}

\newcommand{\rend}[2]{\mathrm{End}({#2}_{#1})}
\newcommand{\lend}[2]{\mathrm{End}({}_{#1}#2)}

\newcommand{\fk}[1]{\mathfrak{#1}}
\newcommand{\Sf}[1]{\mathsf{#1}}
\newcommand{\Bf}[1]{\mathbf{#1}}
\newcommand{\esc}[2]{\langle #1,#2 \rangle}
\newcommand{\Rat}{\mathrm{Rat}}

\newcommand{\T}{\mathcal{T}}

\newcommand{\td}[1]{\widetilde{#1}}
\newcommand{\Sop}[1]{\mathrm{Supp}(#1)}
\newcommand{\scr}[1]{\mathscr{#1}}
\newcommand{\Coind}[1]{{\mathsf{Coind}}^{#1}}
\newcommand{\lr}[1]{\left(\underset{}{} #1 \right)}
\newcommand{\grmod}[1]{\mathsf{gr}\text{-}#1}
\newcommand{\Set}[1]{\left\{\underset{}{} #1 \right\}}

\begin{document}
%\baselineskip 15pt

%Title and Authors
\title{Coinduction functor and simple comodules
\footnote{Research supported by  grant MTM2007-61673 from the Ministerio de
Educaci\'{o}n y Ciencia of Spain, and P06-FQM-01889 from Junta de Andaluc\'{i}a}}
\author{L. El Kaoutit \\
\normalsize Departamento de \'{A}lgebra \\ \normalsize Facultad de Educaci\'{o}n y Humanidades \\
\normalsize  Universidad de Granada \\ \normalsize El Greco N${}^{0}$
10, E-51002 Ceuta, Espa\~{n}a \\ \normalsize
e-mail:\textsf{kaoutit@ugr.es} \and J. G\'omez-Torrecillas \\
\normalsize Departamento de \'{A}lgebra \\ \normalsize Facultad de Ciencias \\
\normalsize Universidad
de Granada\\ \normalsize E18071 Granada, Espa\~{n}a \\
\normalsize e-mail: \textsf{gomezj@ugr.es} }

\date{\today}

\maketitle

\begin{abstract}
Consider a coring with exact rational functor, and a finitely generated and projective right comodule. We construct a functor (\emph{coinduction functor}) which is right adjoint to the hom-functor represented by this comodule. Using the coinduction functor, we establish a bijective map between the set of representative classes of torsion simple right comodules and the set of representative  classes  of simple right modules over the endomorphism ring. A detailed application to a group-graded modules is also given. 
\end{abstract} 

\section*{Introduction}\label{Sec-0}

Let $G$ be a group with neutral element $\Sf{e}$ and $A=\oplus_{x \in G}A_x$ a $G$-graded ring. To each element $x \in G$, one can associate  the restriction functor $(-)_x: \grmod{A} \to \rmod{A_{\Sf{e}}}$ from the category of $G$-graded right $A$-modules to the category of right $A_{\Sf{e}}$-modules. This functor sends an object $ M \in \grmod{A}$ to its homogeneous component $M_x$. Since $(-)_x$ is right exact and commutes with direct sums, a classical result of P. Gabriel tell us that it has a right adjoint functor which we denote by $\Coind{x}:\rmod{A_{\Sf{e}}} \to \grmod{A}$. The functor $\Coind{x}$ is known in the literature as \emph{coinduction functor}, and was first introduced by C. N\u{a}st\u{a}sescu in \cite{Nastasescu:1989a}, for $x =\Sf{e}$. In \cite{Abrams/Menini:1999} G. Abrams and C. Menini defined and studied $\Coind{x}$ in the case of semigroup-graded ring. The use of coinduction functor was crucial to study both simple and injective objects in either a group-graded or semigroup-graded modules categories, see \cite{Nastasescu:1989a, Abrams/Menini:1999, Menini/Nastasescu:1990}.

Recently it was well established \cite{Brzezinski/Wisbauer:2003} that the unified way of studying relative modules and in particular group (semigroup, set-group) graded modules, is the framework of the category of comodules over a suitable coring. The purpose of this paper fits in this direction. Our main aim is to introduce and study the coinduction functor, as well as simple objects, in the context of right comodules over a coring. Some restriction on the base coring are need. Namely, we work with coring which is a member of rational pairing whose associated rational functor is exact. The case of (semi)group-graded modules fits exactly in this situation. A more general case where our methods can perfectly be applied is the case of entwined modules \cite{Brzezinski/Wisbauer:2003}. Of course in this case the factor algebra should be a flat module, while the factor coalgebra should be locally projective module with exact rational functor. Because of space and time, here we only give a detailed application to the case of group-graded modules.

We proceed as follows. In Section \ref{Sec-01}, we give a brief review on rational pairings and rational functors. Section \ref{Sec-1} is devoted to construct the coinduction functors. We consider a coring with exact rational functor. To any finitely generated and projective right comodule, we associate a functor (\emph{coinduction functor}) which is right adjoint to the hom-functor represented by this object (Proposition \ref{AdjointCoinv}). By a classical result of P. Gabriel, we can then deduced that the kernel of this hom-functor is a TTF-class, and the quotient category is equivalent to the category of right modules over the endomorphism ring (Corollary \ref{diag}). Section \ref{Sec-2} presents the main observations on simple comodules. Using the results of Section \ref{Sec-1},  we are able to reconstruct, in a bijective way, any simple right comodule which is torsion with respect to a finitely generated and projective right comodule (Theorem \ref{C-simple}, Corollary \ref{rep-simples}). In the last Section, we give a complete and detailed application to the category of group-graded modules, where we recover some results from \cite{Nastasescu:1989a, Menini/Nastasescu:1990}.

Details on corings and their comodules are easily
available in \cite{Brzezinski/Wisbauer:2003}.

\section{Rational pairing and rational functor}\label{Sec-01}

Rational pairings for coalgebras over commutative rings were
introduced in \cite{Gomez:1998} and used in
\cite{Abuhlail/Gomez/Lobillo:2001} to study the category of right
comodules over the finite dual coalgebra associated to certain
algebras over noetherian commutative rings. This development was
adapted for corings in \cite{ElKaoutit/Gomez/Lobillo:2004c}, see
also \cite{Abuhlail:2003a}. We recall from \cite[Section
2]{ElKaoutit/Gomez/Lobillo:2004c} the definition of this notion: \\
Let $P,Q$ be $A$--bimodules. Any balanced $A$--bilinear form
\begin{displaymath} \esc{-}{-} : P \times Q \longrightarrow A
\end{displaymath} provides in a canonical way two natural transformations $\beta : Q
\tensor{A} - \longrightarrow \hom{}{{}_AP}{-}$ and $\alpha : -
\tensor{A} P \longrightarrow \hom{}{Q_A}{-}$.
Moreover, if $M$ is an $A$--bimodule then $\beta_M$ and $\alpha_M$
are bimodule morphisms. The canonical isomorphisms provide two
bimodule maps
\begin{equation}
\xymatrix@R=0pt{ \beta_A : Q \ar@{->}[r] & \hom{}{_{A}P}{{}_{A}A}
= {}^*P ,\\ q \ar@{->}[r]  & \left[ p \mapsto \esc{p}{q} \right] }
\quad \xymatrix@R=0pt{ \alpha_A : P \ar@{->}[r] &
\hom{}{Q_A}{A_A}= Q^*
\\ p \ar@{->}[r] & \left[ q \mapsto \esc{p}{q} \right] }
\end{equation}
which are bimodule morphisms. So we can recover the balanced
bilinear form if one of the natural transformations is given. The
data $\T = (P,Q,\esc{-}{-})$ are called a \emph{right rational
system over $A$} if $\alpha_M$ is injective for each right
$A$--module $M$, and a \emph{left rational system} if $\beta_N$ is
injective for every left $A$--module $N$. As was mentioned in \cite[Remark 2.2]{ElKaoutit/Gomez/Lobillo:2004c}, 
if $\T=(P,Q,\esc{-}{-})$ is a right rational system, then  ${}_AP$ is a flat module.

\begin{example}\label{locll-proj}
Recall from \cite[Theorem 2.1]{Zimmermann-Huisgen:1976} that a
left $A$--module $P$ is said to be \emph{locally projective} if
for every finite number of elements $p_1,\cdots,p_k\, \in \, P$
there exists a finite system $\{(x_i,\varphi_i)\}_{1 \leq i \leq
n} \subset P \times {}^*P$ such that $$p_j\,\,=\,\,\sum_{i=1}^n
\esc{p_j}{\varphi_i}x_i, \text{ for } 1 \leq j \leq k \text{ (here
} \varphi_i(p_j)=\esc{p_j}{\varphi_i}).$$ For instance, if $P$ is
an $A$--bimodule which is locally projective left $A$--module,
then one can easily show that the three-tuple $(P, {}^*P,
\esc{-}{-})$ where $\esc{-}{-}: P \times {}^*P \rightarrow A$
is the evaluation map, is a right rational system, see
\cite[Lemma 1.29]{Abuhlail:2003a}. 
\end{example}

A \emph{right rational pairing over $A$} is a right rational
system $\T=(\coring{C},B,\esc{-}{-})$ over $A$ consisting of an
$A$--coring $\coring{C}$ and an $A$--ring $B$ (i.e., $B$ is a ring
extension of $A$) such that $\beta_A: B \rightarrow
{}^*\coring{C}$ is a ring anti-morphism where ${}^*\coring{C}$ is
the left dual convolution ring of $\coring{C}$ defined in
\cite[Proposition 3.2]{Sweedler:1975}. As one can easily observe,
the rational pairings are particular instances of right coring
measurings introduced recently in \cite{Brzezinski:2004a}.

\begin{example}\label{par-canonico}
Let $(\coring{C}, \Delta_{\coring{C}}, \varepsilon_{\coring{C}})$
be an $A$--coring such that ${}_A\coring{C}$ is a locally
projective left module. Consider the left colinear endomorphism
ring $\lend{\coring{C}}{\coring{C}}$ as a subring of the linear
endomorphism ring $\lend{A}{\coring{C}}$, that is with
multiplication opposite to the usual composition. Since
$\Delta_{\coring{C}}$ is left $\coring{C}$--colinear and right
$A$--linear map, the canonical ring extension $A \rightarrow
\lend{A}{\coring{C}}$ factors throughout the extension
$\lend{\coring{C}}{\coring{C}} \hookrightarrow
\lend{A}{\coring{C}}$. Therefore, the three-tuple
$\T=(\coring{C},\lend{\coring{C}}{\coring{C}},\esc{-}{-})$ where
the balanced $A$-bilinear $\esc{-}{-}$ is defined by 
$$\esc{c}{f}= \varepsilon_{\coring{C}}(f(c)),\,\, \text{ for every } (c,f) \in \coring{C} \times
\lend{\coring{C}}{\coring{C}}$$ 
is a rational pairing since
$\lend{\coring{C}}{\coring{C}}$ is already a ring anti-isomorphic
to ${}^*\coring{C}$ via the beta map associated to $\esc{-}{-}$.
We refer to $\T$ as \emph{the right canonical pairing} associated
to $\coring{C}$.
\end{example}

Given $\T=(\coring{C},B,\esc{-}{-})$ any right rational pairing
over $A$,  one can define a functor called a \emph{right rational
functor} as follows. Let $M$ be any right $B$--module. An element
$m \in M$ is called \emph{rational} if there exists a set of
\emph{right rational parameters} $\{(c_i,m_i)\} \subseteq
\coring{C} \times M$ such that $ m b = \sum_i m_i\esc{c_i}{b}$,
for all $b \in B$. The set of all rational elements in $M$ is
denoted by $\Rat^{\T}(M)$. As it was explained in \cite[Section
2]{ElKaoutit/Gomez/Lobillo:2004c}, the proofs detailed in
\cite[Section 2]{Gomez:1998} can be adapted in a straightforward
way in order to get that $\Rat^{\T}(M)$ is a $B$--submodule of $M$
and the assignment $M \mapsto \Rat^{\T}(M)$ is a well defined
functor
\begin{equation*}
\Rat^{\T} : \rmod{B} \rightarrow \rmod{B},
\end{equation*}
which is in fact a left exact preradical. Therefore, the full
subcategory $\Rat^{\T}(\rmod{B})$ of $\rmod{B}$ whose objects are
those $B$--modules $M$ such that $\Rat^{\T}(M) = M$ is a closed
subcategory. Furthermore, $\Rat^{\T}(\rmod{B})$ is a Grothendieck
category which is shown to be isomorphic to the category of right
comodules $\rcomod{\coring{C}}$ as \cite[Theorem
2.6']{ElKaoutit/Gomez/Lobillo:2004c} asserts. In this way, we still denote $$\Rat^{\T}: \rmod{B} \longrightarrow \Rat^{\T}(\rmod{B}) \cong \rcomod{\coring{C}}.$$

\section{Coinduction functor in corings}\label{Sec-1}

Let $\T\,=\,(\coring{C}, B,\esc{-}{-})$ be a right rational
pairing over $A$, and $\Rat^{\T}: \rmod{B} \rightarrow
\rcomod{\coring{C}}$ the associated rational functor. We know that
there is an adjunction $$\xymatrix{ \Rat^{\T}: \rmod{B}
\ar@<0,5ex>[r] & \ar@<0,5ex>[l] \rcomod{\coring{C}}: i^{\T},  }$$
where $i^{\T}$ is left adjoint to $\Rat^{\T}$. Consider
$(\Sigma,\rho_{\Sigma})$ a right $\coring{C}$-comodule such that
$\Sigma_A$ is finitely generated and projective module. Let us
denote by $T_{\Sigma}=\rend{\coring{C}}{\Sigma}$ its colinear
endomorphism ring. As we have seen in Section \ref{Sec-01},  we also have $T_{\Sigma}=\rend{B}{\Sigma}$. It
is well known that the coinvariant functor
$\hom{\coring{C}}{\Sigma}{-}: \rcomod{\coring{C}} \rightarrow
\rmod{T_{\Sigma}}$ is a right adjoint to the tensor product
functor $-\tensor{T_{\Sigma}}\Sigma: \rmod{T_{\Sigma}}
\longrightarrow \rcomod{\coring{C}}$. We are interested in looking
at a possible right adjoint of the functor
$\hom{\coring{C}}{\Sigma}{-}$. So, if this right adjoint functor
exists, then $\Sigma$ should be a projective right
$\coring{C}$-comodule, since $\rcomod{\coring{C}}$ is an abelian
category (recall that ${}_A\coring{C}$ is flat by the pairing
$\T$). Henceforth, we assume that $\Sigma_{\coring{C}}$ is a
finitely generated and projective comodule. Since $i^{\T}$ is left
adjoint to $\Rat^{\T}$, if $\Rat^{\T}$ is an exact functor, then
$i^{\T}(\Sigma_{\coring{C}}) =\Sigma_{B}$ is finitely generated
and projective right $B$-module.  We denote by
$\Sigma^{\star}=\hom{B}{\Sigma}{B}$ its right dual. Then by the rational pairing $\T$, the structure
of $(T_{\Sigma},\coring{C})$-bicomodule on $\Sigma$ is equivalent to the structure of
$(T_{\Sigma},B)$-bimodule, and so $\Sigma^{\star}$ becomes a
$(B,T_{\Sigma})$-bimodule.

\begin{proposition}\label{AdjointCoinv}
Let $\T=(\coring{C},B,\esc{-}{-})$ be a right rational pairing over a ring $A$, and $(\Sigma,\rho_{\Sigma})$
a finitely generated and projective right $\coring{C}$-comodule with endomorphism ring $T_{\Sigma}$. If the rational functor $\Rat^{\T}: \rmod{B} \rightarrow
\rcomod{\coring{C}}$ is exact, then $\Sigma_{B}$ is finitely generated and projective module. In particular, the functor $$\Rat^{\T}
\,\circ\, \hom{T_{\Sigma}}{\Sigma^{\star}}{-}: \rmod{T_{\Sigma}}
\longrightarrow \rcomod{\coring{C}}$$ is right adjoint to the
coinvariant functor $\hom{\coring{C}}{\Sigma}{-}:
\rcomod{\coring{C}} \rightarrow \rmod{T_{\Sigma}}$.
\end{proposition}
\begin{proof}
We have seen in the preamble of this Section that $\Sigma_{B}$ should be finitely generated and projective. Therefore, the stated adjunction follows from the following composition of adjoint pairs of
functors $$\xymatrix@C=80pt{ \rmod{T_{\Sigma}} \ar@<0,5ex>[r]^{
\hom{T_{\Sigma}}{\Sigma^{\star}}{-} } &
\ar@<0,5ex>[l]^{\hom{B}{\Sigma}{-}} \rmod{B}
\ar@<0,5ex>[r]^{\Rat^{\T}} & \ar@<0,5ex>[l]^{i^{\T}}
\rcomod{\coring{C}}, }$$ since we know that $
\hom{B}{\Sigma}{-} \circ i^{\T}   \,= \, \hom{\coring{C}}{\Sigma}{-}$.
\end{proof}

From now on, we denote by 
$$\Coind{\Sigma}(-):= \Rat^{\T} \,\circ\,
\hom{T_{\Sigma}}{\Sigma^{\star}}{-}: \rmod{T_{\Sigma}} \longrightarrow \rcomod{\coring{C}}$$ 
and refer to  as \emph{the
coinduction functor} associated to the finitely generated and
projective comodule $\Sigma_{\coring{C}}$. Fix a right dual basis
$\{(u_j,\,u_j^{\star})\}_j$ for $\Sigma_B$. The unit and counit of
the adjunction stated in Proposition \ref{AdjointCoinv} are given as follows:
\begin{equation}\label{unit-counit}
\xi_{Y_{T_{\Sigma}}}: \hom{\coring{C}}{\Sigma}{\Coind{\Sigma}(Y)}
\longrightarrow Y, \quad \left( f \longmapsto
\sum_{j}f(u_j)(u_j^{\star}) \right),$$
$$ \eta_{X_{\coring{C}}}: X \longrightarrow
\Coind{\Sigma}\left(\hom{\coring{C}}{\Sigma}{X}\right), \quad
\left( x \longmapsto \left[ \underset{}{} u^{\star} \mapsto[v
\mapsto xu^{\star}(v)] \right]
 \right)
\end{equation}
for every comodule $X_{\coring{C}}$ and module $Y_{T_{\Sigma}}$.
The counit $\xi_{-}$ is actually a natural isomorphism. Effectively, for
every right $T_{\Sigma}$-module $Y$, we have
\begin{eqnarray*}
\hom{\coring{C}}{\Sigma}{\Coind{\Sigma}(Y)} &=& \hom{B}{\Sigma}{\hom{T_{\Sigma}}{\Sigma^{\star}}{Y}} \\
&\cong & \hom{T_{\Sigma}}{\Sigma\tensor{B}\Sigma^{\star}}{Y}, \quad T_{\Sigma} \cong
\Sigma\tensor{B}\Sigma^{\star} \\
& \cong & \hom{T_{\Sigma}}{T_{\Sigma}}{Y} \,\, \cong \,\,
Y_{T_{\Sigma}},
\end{eqnarray*}
and the composition of those isomorphisms is exactly $\xi_{Y}$ defined in \eqref{unit-counit}.

Let us denote 
$$ \scr{C}_{\Sigma}\, =\,
\mathrm{Ker}(\hom{\coring{C}}{\Sigma}{-})\, =\, \left\{ L \in
\rcomod{\coring{C}} |\, \hom{\coring{C}}{\Sigma}{L} = 0\right\}.$$
Consider 
$$\scr{T}_{\Sigma} \,=\, \left\{ M \in
\rcomod{\coring{C}} |\, \hom{\coring{C}}{M}{L} = 0, \,\forall L
\in \, \scr{C}_{\Sigma}\right\}$$ the torsion class associated to the torsion-free class $\scr{C}_{\Sigma}$. The corresponding
idempotent radical is  
\begin{equation}\label{radical}
\fk{r}_{\Sigma}: \rcomod{\coring{C}}
\longrightarrow \rcomod{\coring{C}}, \,\,\, \left( M \mapsto
\fk{r}_{\Sigma}(M)\, =\,\sum \{ \mathrm{Im}(f)|\, f: N
\rightarrow M,\, N \in \, \scr{T}_{\Sigma} \} \right).
\end{equation}

The following corollary is a direct consequence of
\cite[Proposition 5, p. 374]{Gabriel:1962}, since
$\Sigma_{\coring{C}}$ is projective and $\xi_{-}$ is a natural
isomorphism. For the definition of TTF-class in Grothendieck
categories we refer the reader to \cite[Chap. VI]{Stenstrom:1975}.

\begin{corollary}\label{diag}
With assumptions as in Proposition \ref{AdjointCoinv},  the full subcategory $\scr{C}_{\Sigma}$ is a TTF-class, and there
is a commutative diagram
$$\xymatrix@R=40pt@C=60pt{ \rcomod{\coring{C}}
\ar@{->}^-{\Bf{T}_{\Sigma}}[r]
\ar@{->}_-{\hom{\coring{C}}{\Sigma}{-}}[dr] &
\rcomod{\coring{C}}/\scr{C}_{\Sigma}
\ar@{->}^-{\Bf{R}_{\Sigma}}[d]  \\ & \rmod{T_{\Sigma}}, }$$ where
$\Bf{T}_{\Sigma}$ is the localizing functor and $\Bf{R}_{\Sigma}$
is an equivalence of categories with inverse
$\Bf{L}_{\Sigma}=\Bf{T}_{\Sigma} \circ \Coind{\Sigma}$.
\end{corollary}

\section{Simple comodules}\label{Sec-2}
In this section we make some observations on simple  right
$\coring{C}$-comodules, which are also $\scr{T}_{\Sigma}$-torsion comodules.
To this end the following two lemmata will be needed.

\begin{lemma}\label{loc-simple}
Let $\cat{G}$ be a Grothendieck category and $\scr{C}$ a TTF-class
of $\cat{G}$ with torsion theory $(\scr{T},\scr{C})$ and
associated radical functor $\fk{r}: \cat{G} \rightarrow \cat{G}$.
Let $\xymatrix{ \Bf{T} : \cat{G}\ar@<0,5ex>[r] & \ar@<0,75ex>[l]
\cat{G}/\scr{C}: \Bf{S} }$ be the canonical adjunction of
localization. Consider an object $N$ of $\cat{G}$ such that
$\fk{r}(N) \neq 0$. If $\Bf{T}(N)$ is a simple object of
$\cat{G}/\scr{C}$, then $\fk{r}(N)$ is a simple object of
$\cat{G}$ as well.
\end{lemma}
\begin{proof}
Let $\alpha: M\hookrightarrow \fk{r}(N)$ be a non zero
monomorphism of $\cat{G}$. Since $\Bf{T}$ is an exact functor, we then
get a monomorphism $\Bf{T}(\alpha): \Bf{T}(M) \hookrightarrow
\Bf{T}(\fk{r}(N))$. Clearly $\Bf{T}(\alpha) \neq 0$ since
$\mathrm{Im}(\alpha)$ is not an object of $\scr{C}$. On the other hand, we know that
$N/\fk{r}(N) \in \scr{C}$, which means that  $\Bf{T}(\fk{r}(N)) \cong
\Bf{T}(N)$. Hence $\Bf{T}(\alpha)$ is an isomorphism, since
${\bf T}(N)$ is assumed to be simple. Consequently $\mathrm{Coker}(\alpha) = \fk{r}(N) / \mathrm{Im}(\alpha) \in \scr{C}$ which
implies that $N/\mathrm{Im}(\alpha) \in \scr{C}$. Whence
$\fk{r}(N)=\mathrm{Im}(\alpha)$, and so $\alpha$ is an isomorphism.
Therefore, $\fk{r}(N)$ is a simple object of $\cat{G}$.
\end{proof}

\begin{lemma}\label{mod-simple}
Let $B$ be any ring and $\Sigma_B$ a finitely generated and
projective module with endomorphism ring $T=\rend{B}{\Sigma}$ and
right dual module $\Sigma^{\star}=\hom{B}{\Sigma}{B}$. Consider a maximal
right ideal $I_B$ of $B_B$ and assume that $(I\Sigma^{\star})_T
\lneqq \Sigma^{\star}_T$. Then $(I\Sigma^{\star})_T$ is a maximal
submodule of $\Sigma^{\star}_T$.
\end{lemma}
\begin{proof}
It is sufficient to show that, for any element $u^{\star} \in \Sigma^{\star} \setminus I\Sigma^{\star}$, we
have 
\begin{equation}\label{max}
I\Sigma^{\star} \,\, + \,\, u^{\star} \, T \,\, = \,\,
\Sigma^{\star}.
\end{equation}

Let $\{u_i,u_i^{\star}\}_{1 \leq i \leq n}$ be a finite dual basis
for $\Sigma_B$, and choose an arbitrary element $u^{\star} \in\Sigma^{\star} \setminus I\Sigma^{\star}$. We know that $u^{\star} =
\sum_{1}^{n}u^{\star}(u_i)u_i^{\star}$ and $\{1,\cdots,n\}=
\Lambda \uplus \Lambda'$ where $(i \in \Lambda \, \Rightarrow\,
u^{\star}(u_i) \in I)$ and $ (i \in \Lambda' \, \Rightarrow \,
u^{\star}(u_i) \notin I)$ (by assumption we have $\Lambda
\varsubsetneqq \{1,\cdots,n\}$). Since $I_B$ is maximal, we have
$$I \,\, + \,\, \sum_{i \in \,\Lambda'} u^{\star}(u_i) \, B \,\, =
\,\, B.$$ Hence
\begin{equation}\label{max'}
I\Sigma^{\star} \,\, + \,\, \sum_{i \in \,\Lambda'}
u^{\star}(u_i)\Sigma^{\star} \,\, = \,\, \Sigma^{\star}\quad
\text{ in the category of modules } \rmod{T}.
\end{equation}
Consider $K:=\sum_{i \in \,\Lambda'} u^{\star}(u_i)\Sigma^{\star}$
as right $T$-module. It is easily checked that the following
$T$-linear map is a monomorphism
$$\varsigma:\,\, K \longrightarrow u^{\star} \, T,\qquad \left(
\sum_{i \in \,\Lambda'} u^{\star}(u_i)v^{\star}_i \longmapsto
u^{\star}(\sum_{i \in \,\Lambda'}t_i) \right), $$ where the
$t_i$'s are elements of $T$ given by  $$t_i:   \Sigma \longrightarrow \Sigma, \quad [x \mapsto u_iv^{\star}_i(x)], \,\, \forall i\in
\,\Lambda'.$$ Therefore, equality \eqref{max'} implies equality
\eqref{max}.
\end{proof}

\begin{theorem}\label{C-simple}
Let $\T=(\coring{C}, B, \esc{-}{-})$ be a right rational pairing over $A$ with exact rational functor. 
Consider a finitely generated and projective right $\coring{C}$-comodule $\Sigma$ whose endomorphism ring is $T_{\Sigma}$, and 
with associated idempotent radical $\fk{r}_{\Sigma}$ as in
Section \ref{Sec-1}. 
\begin{enumerate}[(a)]
\item If $Y_{T_{\Sigma}}$ is a simple module, then
$\fk{r}_{\Sigma}(\Coind{\Sigma}(Y))_{\coring{C}}$ is a simple
comodule, and 
$$ Y_{T_{\Sigma}} \, \cong\, \hom{\coring{C}}{\Sigma}{\fk{r}_{\Sigma}(\Coind{\Sigma}(Y))}_{T_{\Sigma}}.$$ 
Moreover, given two simple right $T_{\Sigma}$-module $Y$ and $Y'$, we have 
$$ Y_{T_{\Sigma}} \, \cong Y'_{T_{\Sigma}}\, \, \Longleftrightarrow \, \fk{r}_{\Sigma}(\Coind{\Sigma}(Y))_{\coring{C}}\, \cong  \fk{r}_{\Sigma}(\Coind{\Sigma}(Y'))_{\coring{C}}.$$

\item  If $S_{\coring{C}}$ is a simple comodule such
that $\hom{\coring{C}}{\Sigma}{S} \neq 0$. Then
$\hom{\coring{C}}{\Sigma}{S}_{T_{\Sigma}}$ is a simple module.
Furthermore, there is an isomorphism of right comodules
$$ S_{\coring{C}} \,\, \cong\,\,
\fk{r}_{\Sigma}\left(\underset{}{}
\Coind{\Sigma}\left(\underset{}{}\hom{\coring{C}}{\Sigma}{S}\right)\right)_{\coring{C}}.$$
\end{enumerate}
\end{theorem}
\begin{proof}
$(a)$. Suppose that $\fk{r}_{\Sigma}(\Coind{\Sigma}(Y)) = 0$, that is,
$\Coind{\Sigma}(Y) \in \scr{C}_{\Sigma}$. So $$Y_{T_{\Sigma}} \,\,
\cong \,\,\hom{\coring{C}}{\Sigma}{\Coind{\Sigma}(Y)}_{T_{\Sigma}}
=0,\quad (\text{ via } \xi_{Y}, \text{ equation } \eqref{unit-counit}).$$ This is a contradiction since we know that $Y\neq 0$.
Henceforth, $\fk{r}_{\Sigma}(\Coind{\Sigma}(Y)) \neq 0$, and clearly $Y \cong \hom{\coring{C}}{\Sigma}{\fk{r}_{\Sigma}(\Coind{\Sigma}(Y))}$. In order to check that $\fk{r}_{\Sigma}(\Coind{\Sigma}(Y))$ is simple right ${T_{\Sigma}}$-module, it suffices to show by Lemma
\ref{loc-simple} that $\Bf{T}_{\Sigma}(\Coind{\Sigma}(Y))$  is a simple object in the quotient category
$\rcomod{\coring{C}}/\scr{C}_{\Sigma}$. This is fulfilled since by Corollary
\ref{diag}, ${\bf L}_{\Sigma} =\Bf{T}_{\Sigma}\circ\Coind{\Sigma}$ is an equivalence
of categories $\rmod{T_{\Sigma}}$ and
$\rcomod{\coring{C}}/\scr{C}_{\Sigma}$. 
The implication $(\Rightarrow)$ of the stated equivalence is obvious. Let us check the converse one. So assume we are given $Y$ and $Y'$ two simple right $T_{\Sigma}$-modules such that $\fk{r}_{\Sigma}(\Coind{\Sigma}(Y)) \cong  \fk{r}_{\Sigma}(\Coind{\Sigma}(Y'))$. Applying the functor ${\bf T}_{\Sigma}$ of Corollary \ref{diag}, we obtain the following chain of isomorphisms 
$${\bf T}_{\Sigma}\lr{\Coind{\Sigma}(Y)} \cong {\bf T}_{\Sigma}\lr{\fk{r}_{\Sigma}(\Coind{\Sigma}(Y))} \cong  {\bf T}_{\Sigma}\lr{\fk{r}_{\Sigma}(\Coind{\Sigma}(Y'))} \cong {\bf T}_{\Sigma}\lr{\Coind{\Sigma}(Y')},$$ where we have used the same argument as in the proof of Lemma \ref{loc-simple} for the first and the third isomorphisms. Now, since ${\bf T}_{\Sigma}\circ \Coind{\Sigma}$ is an equivalence of categories, we then get $Y \cong Y'$ as right $T_{\Sigma}$-modules. 

$(b)$. Since $\Sigma_B$ is finitely generated and projective module by Proposition \ref{AdjointCoinv}, we
have $\hom{\coring{C}}{\Sigma}{S}_{T_{\Sigma}}$ is simple module if and only if
$(S\tensor{B}\Sigma^{\star})_{T_{\Sigma}}$ it is. On the other
hand, $S_B\cong (B/I)_B$, for some maximal right ideal $I_B$ of
$B_B$. Hence,
$$ (S\tensor{B}\Sigma^{\star})_{T_{\Sigma}} \,\,
\cong \left(\underset{}{} (B/I) \tensor{B} \Sigma^{\star}
\right)_{T_{\Sigma}} \,\, \cong \,\, \left(\underset{}{}
\Sigma^{\star}/I\Sigma^{\star} \right)_{T_{\Sigma}},$$ and so
$\hom{B}{\Sigma}{S} \neq 0$ implies
$(I\Sigma^{\star})_{T_{\Sigma}} \lneqq
\Sigma^{\star}_{T_{\Sigma}}$. The hypothesis of Lemma
\ref{mod-simple} is then fulfilled, which implies that
$(\Sigma^{\star}/I\Sigma^{\star})_{T_{\Sigma}}$ is simple module
and so is $(S\tensor{B}\Sigma^{\star})_{T_{\Sigma}}$. Therefore,
$\fk{r}_{\Sigma}\left(\underset{}{}\Coind{\Sigma}
\left(\underset{}{}\hom{\coring{C}}{\Sigma}{S}\right)\right)$ is
simple right $\coring{C}$-comodule, by item (a).
Lastly, the stated isomorphism of comodules is deduced from the natural transformation $\eta_{-}$ of equation \eqref{unit-counit} as follows
$$\fk{r}_{\Sigma}(\eta_{S}): S=\fk{r}_{\Sigma}(S) \,\, \cong \,\,
\fk{r}_{\Sigma}\left(\underset{}{}\Coind{\Sigma}
\left(\underset{}{}\hom{\coring{C}}{\Sigma}{S}\right)\right).$$
\end{proof}

Let us consider the following two sets. $\scr{S}_{\Sigma}$ the set of isomorphism classes $[S]$ represented by simple right $\coring{C}$-comodule $S$ such that $\hom{\coring{C}}{\Sigma}{S} \neq 0$. While $\scr{S}_{T_{\Sigma}}$ is the set of isomorphism classes $[Y]$ represented by simple right $T_{\Sigma}$-modules $Y$. Using Theorem \ref{C-simple}, we then get 

\begin{corollary}\label{rep-simples}
Let $\T=(\coring{C}, B, \esc{-}{-})$ be a right rational pairing over $A$ with exact rational functor. 
Consider a finitely generated and projective right $\coring{C}$-comodule $\Sigma$ whose endomorphism ring is $T_{\Sigma}$. Then there is a bijective map given by
$$
\xymatrix@R=0pt{ \scr{S}_{\Sigma} \ar@{->}[rr] & & \scr{S}_{T_{\Sigma}} \\ [S] \ar@{|->}[rr] & & \left[ \hom{\coring{C}}{\Sigma}{S}\right] \\ \left[ \fk{r}_{\Sigma} \Coind{\Sigma}(Y)\right]  & &  \ar@{|->}[ll] [Y]. }
$$ 
\end{corollary}

\section{Application to Group-graded modules}
We show that the coinduction functor 
\cite{Nastasescu:1989a} in group-graded modules involves a
rational functor associated to a suitable rational pairing. This
pairing is constructed using the smash products
\cite{Quinn:1985,Cohen/Montgomery:1984} and the canonical coring
arising from a group-graded base ring. Using this functor, we will apply the results of Section \ref{Sec-2} to simple group-graded right modules.

In what follows, we consider a group $G$ with neutral element
$\Sf{e}$ and a $G$-graded base ring $A=\oplus_{x \in \,G}A_x$. It
is well known that the free left $A$-module with basis $G$,
denoted by $\coring{C}=AG$, admits a right $A$-action given by the
rule $xa_y\,=\, a_y(xy)$ for every homogeneous element $a_y \in
A_y$ and every $x,y \in G$. It turns out that $\coring{C}$ is in
fact an $A$-bimodule equipped with an $A$-coring structure with
comultiplication $\Delta(ax)= ax \tensor{A}x$ and counit
$\varepsilon(ax) =a$, for every $x \in G$ and $a \in A$. Of course
$G$ is contained in the set of all grouplike elements of
$\coring{C}$; recall that a grouplike element of an $A$-coring
$\coring{C}$ is an element $g \in \coring{C}$ such that
$\varepsilon_{\coring{C}}(g)=1$ and
$\Delta_{\coring{C}}(g)=g\tensor{A}g$.

In this way the category of all right graded $A$-module
$\grmod{A}$ is isomorphic to the category of all right
comodules $\rcomod{\coring{C}}$. This isomorphism, over objects, is
given as follows. For every $x \in G$, we consider the right
$\coring{C}$-comodule $([x]A,\rho_{[x]A})$ whose underlying right
$A$-module coincides with $A_A$, i.e. $([x]A)_A=A_A$ and its
coaction map is  $\rho_{[x]A}: [x]A \to  \coring{C} \cong
A\tensor{A}\coring{C}$ sending  $a \mapsto x a $. Given a right comodule
$(M,\rho_M)$, then $M$ decomposes as $M=\oplus_{x \in
\,G}M^{cov(x)}$, where $$M^{cov(x)}=\{m \in M|\,
\rho_{M}(m)=m\tensor{A}x\}, \text{ for every }x \in G.$$ One can easily
check that this gives to $M$ a structure of right graded
$A$-module, and there is an isomorphism of abelian groups
\begin{equation}\label{coinv-x}
M^{cov(x)} \cong \hom{\coring{C}}{[x]A}{M},\,\, \text{ for every } x \in\, G.
\end{equation}
In particular, we have 
$$ \hom{\coring{C}}{[x]A}{[y]A}\,=\, \Set{a \in A|\, ax=ya}\,=\, A_{xy^{-1}},\,\, \text{ for every } x, y \in\, G,$$ which implies that the endomorphism ring of each right comodules $[x]A$ is isomorphic to $A_{\Sf{e}}$. That is, $${\rm End}_{\coring{C}}([x]A)\,\cong \, A_{\Sf{e}},\,\, \text{ for every } x \in\, G.$$  

Conversely, given $N=\oplus_{x \in \,G}N_x$ a graded right
$A$-module, one can endow this module by a right $\coring{C}$-coaction
sending $\rho_N(n_x)=n_x \tensor{A}x$, for every homogeneous
element $n_x \in N_x$, $x \in G$. The isomorphism of categories $\rcomod{\coring{C}} \cong \grmod{A}$ is now established taking into account that its action on morphisms is an identity. Moreover, for every element $x \in G$, we have a commutative
diagram of functor
$$\xymatrix@C=100pt@R=30pt{ \rmod{A_{\Sf{e}}}
 & \rcomod{\coring{C}}  \ar@{->}_-{\hom{\coring{C}}{[x]A}{-}}[l] \\
& \ar@{->}^-{(-)_{x}}[lu] \ar@{->}_-{\cong}[u]  \grmod{A}, }$$ where
$(-)_{x}: \grmod{A} \rightarrow \rmod{A_{\Sf{e}}}$ is the
induction functor which sends 
every right graded $A$-module $N=\oplus_{x \in \,G}N_x$ to $(N)_{x} = N_{x}$ its homogeneous component of degree $x$.

On the other hand, for each $x \in G$, the left $A$-submodule
$P_x=Ax$ of $\coring{C}$ generated by $x$ is clearly a left
$\coring{C}$-comodule with coaction sending $\lambda_{P_x}(ax)=ax
\tensor{A} x \in \coring{C} \tensor{A} P_x$, $a \in  A$.
Therefore, ${}_{\coring{C}}\coring{C}= \oplus_{x \in \,G}P_x$ is a
direct sum of left $\coring{C}$-comodules. For any $x \in G$, we
denote by $e_x: \coring{C} \twoheadrightarrow P_x \hookrightarrow
\coring{C}$ the composition of the canonical injection and
projection in this direct sum. The $e_x$'s form a set of
orthogonal idempotent elements in the endomorphism ring
$\lend{\coring{C}}{\coring{C}}$. The ring extension
$\widetilde{(-)}: \,A \rightarrow \lend{\coring{C}}{\coring{C}}$ which sends any element $a \in A$ to the map $\widetilde{a}: c \mapsto ca$, is a monomorphism as ${}_A\coring{C}$ is a free
module. We thus identify $A$ with its image $\widetilde{A}$.

\begin{proposition}\label{smashProduct}
Let $A$ be a $G$-graded ring by a group $G$. Consider
the associated $A$-coring $\coring{C}=AG$, and denote by $B$ the
sub-ring of $\lend{\coring{C}}{\coring{C}}$ generated by
$\widetilde{A}$ and the set of orthogonal idempotents $\{e_x\}_{x
\in\,G}$. Then
\begin{enumerate}[(i)]
\item $B_A$ is a free right module with basis $\{1,e_x\}_{x \in \,
G}$, that is, $B = \td{A} \oplus \left( \oplus_{x \in
\,G}\,e_x\td{A}\right)$.

\item There is a right rational pairing
$\T=(\coring{C},B,\esc{-}{-})$, where the bilinear form is defined
by the rule:  $$\esc{x}{e_y} \,=\, \delta_{x,\,y} \text{ (the
Kronecker delta) and } \,\, \esc{x}{\td{a}} \,=\, a,$$ for every
$x,\, y \in \,G$, $a \in A$, and satisfies
$$c \,\, = \,\,  \sum_{z \in \, G} \esc{c}{e_z} \, z, \text{ (finite
sum) for every } c \in \coring{C}.$$

\item Let $\Rat^{\T}:\rmod{B} \rightarrow \rcomod{\coring{C}}$ be
the rational functor associated to the pairing $\T$ of item
$(ii)$. Then the trace ideal is $\Rat^{\T}(B_B) \, =\, \oplus_{x
\in \,G}\,e_x \td{A}$, and the functor $\Rat^{\T}$ is exact. In particular, for every right $B$-module $M$, we have
$$\Rat^{\T}(M_B) \,\,= \,\, M\,\Rat^{\T}(B_B) \,\, = \,\, \oplus_{x \in \, G}
\,Me_x.$$
\end{enumerate}
\end{proposition}
\begin{proof}
$(i)$. This was proved in \cite{Quinn:1985} and
\cite{Nastasescu:1989a}, using the fact that $\lend{A}{\coring{C}}$ is
isomorphic to the ring of finite rows $|G| \times |G|$-matrices
over $A$. Here, we can restrict our selves to the colinear
endomorphism ring $\lend{\coring{C}}{\coring{C}}$. The statement of this item follows then from the following two facts. The first is that, for every $x,y \in \,G$, $a \in A$, $a_y \in A_y$, we have
\begin{equation}\label{a-e}
  e_x \,\, \td{a} \,\, e_y \,\,=\,\, \td{\pi_{y^{-1}x}(a)} \,
e_y,\qquad e_x \, \td{a_y}  \,\,=\,\, \td{a_y} \, e_{xy},
\end{equation}
where $\pi_z: A\rightarrow A_z$, $z \in \, G$ are the canonical projections. The second is that, if $a \in A$ and  $e_x \td{a}= 0$, for some $x \in G$, then $a=0$.

$(ii)$. To show that $\esc{-}{-}$ is well defined, we only prove
that $\esc{-}{-}$ is $A$-balanced. So, let $a = \sum_{z \in
\,G}\pi_z(a) \in A$ and $x,y \in G$, we have
$$ \esc{xa}{e_y} \,\, = \,\, \sum_{z \in \, G} \pi_z(a)
\esc{xz}{e_y} \,\, =\,\, \sum_{z \in \,G}\pi_z(a) \delta_{xz,\,y}
\,\, = \,\, \pi_{x^{-1}y}(a),$$ and $$ \esc{x}{\td{a}e_y} \,\, =
\,\, \sum_{z \in \, G} \esc{x}{\td{\pi_z(a)}e_y} \,\, =\,\,
\sum_{z \in \,G} \esc{x}{e_{yz^{-1}}} \pi_z(a)  \,\, = \,\,
\sum_{z \in \,G}\pi_z(a) \delta_{x,\, yz^{-1}} \,\,=\,\,
\pi_{x^{-1}y}(a),$$ where we have used equation \eqref{a-e}.
Hence, $\esc{xa}{e_y} = \esc{x}{\td{a}e_y}$, for every $x,y \in
\,G$ and $a \in A$ which implies that $\esc{-}{-}$ is well
defined. Given $c \in \,\coring{C}$, we denote as usual by
$\Sop{c} \subseteq G$ (finite subset) the support of $c$, that is
$z \in \Sop{c} \Leftrightarrow \kappa_z(c) \neq 0$, where
$\kappa_x: \coring{C} \rightarrow Ax$, $x \in \,G$ are the
canonical projections. Obviously, $c =\sum_{z \in \, \Sop{c}}
a^{c,\,z} z$ with $a^{c,\,z} \in A$, and if $z_0 \in \Sop{c}$,
then $$\esc{c}{e_{z_0}}= \sum_{z \in \Sop{c}}
a^{c,\,z}\delta_{z,\,z_0} = a^{c,\,z_0}.$$ Therefore, $c = \sum_{z
\in\, \Sop{c}}\esc{c}{e_z} z$ as claimed. Now, the natural
transformation associated to this bilinear form is then given by
$$ \alpha_M: M\tensor{A} \coring{C} \longrightarrow
\hom{A}{B}{M},\,\, (m \tensor{A}c \longmapsto [b \mapsto
m\esc{c}{b}]), $$ for every right $A$-module $M$. Suppose that
$\alpha_M(\sum_i^nm_i\tensor{A}c_i) =0$, and put $H=
\cup_{i=1}^{n} \mathrm{Supp}(c_i)$, then $c_i = \sum_{z \in \,
H}\esc{c_i}{e_z} z$, for all $i=1,\cdots,n$. We thus obtain
$$\sum_{i=1}^n m_i \tensor{A} c_i \,\, = \,\, \sum_{i=1}^n
\sum_{z \in \,H}m_i \tensor{A} \esc{c_i}{e_z}z \,\, = \,\, \sum_{z
\in\,H}\left( \sum_{i=1}^n m_i \esc{c_i}{e_z}\right)\tensor{A} z
\,\, =\,\, 0,$$ since $\alpha_M(\sum_i^nm_i\tensor{A}c_i)(e_z) =
\sum_im_i\esc{c_i}{e_z} =0$, for all $z \in \,G$. This shows that
$\alpha_{M}$ is injective for every right $A$-module $M$. The
remainder of axioms are obviously satisfied, which shows that $\T$ is rational
pairing.

$(iii)$. By equations \eqref{a-e}, we know that $\oplus_{x \in \,G}e_x\td{A}$
is a two-sided idempotent ideal of $B$ which is left and right
pure $B$-submodule of $B$. Now, for every $x \in \,G$,
$(e_x\td{A},\rho_{e_x\td{A}})$ is a right $\coring{C}$-comodule
with coaction $\rho_{e_x\td{A}}: e_x\td{A} \rightarrow e_x\td{A}
\tensor{A} \coring{C}$ sending $e_x\td{a} \mapsto e_x \tensor{A}
xa$. Thus, we have an inclusion $\oplus_{x \in \, G}e_x\td{A}
\,\subseteq \Rat^{\T}(B_B)$. Using the definition of rational elements Section \ref{Sec-01}, item $(i)$ and the second equation of \eqref{a-e}, we can show that $ \Rat^{\T}(B_B) \cap \td{A} = \{0\}$. Therefore, $\oplus_{x \in \, G}e_x\td{A} \,= \, \Rat^{\T}(B_B)$. The
rational functor $\Rat^{\T}$ is now exact by \cite[Theorem 1.2]{ElElKaoutit/Gomez:2008}  as $\coring{C}\, \Rat^{\T}(B_B) =\coring{C}$. Using
again \cite[Theorem 1.2]{ElElKaoutit/Gomez:2008} and equations \eqref{a-e}, we then obtain 
$\Rat^{\T}(M)=M \,\Rat^{\T}(B_B) = \oplus_{x \in \,G} Me_x$, for
every right $B$-module $M$.
\end{proof}

\begin{remark}
By \cite[Theorem 1.2]{ElElKaoutit/Gomez:2008} and  Proposition \ref{smashProduct}, we know that
the category $$\Rat^{\T}(\rmod{B})=\{ M \in \rmod{B}|\,\,
M=\oplus_{x \in \,G}Me_x\}$$ is a localizing sub-category of
$\rmod{B}$, this is \cite[Proposition 1.1]{Nastasescu:1989b}.
The adjoint pair of functors constructed in \cite[Section
1]{Nastasescu:1989b}, can be 
easily obtained as follows. From the proof
of Proposition \ref{smashProduct}, it is easily seen that the
right $A$-modules $e_x\td{A}$ admit a structure of right
$B$-module coming from its right $\coring{C}$-coaction. Moreover,
$e_x\td{A}$ and  $[x]A$ are isomorphic as right
$\coring{C}$-comodules, for every $x \in \,G$. Thus the trace
ideal $\fk{a}=\Rat^{\T}(B_B) \cong \oplus_{x \in \,G}[x]A$ as
right $\coring{C}$-comodule. Therefore, the right adjoint functor
$F = \hom{B}{\fk{a}}{-} \circ i^{\T}: \rcomod{\coring{C}}
\rightarrow \rmod{B} $ of $\Rat^{\T}: \rmod{B} \rightarrow
\rcomod{\coring{C}}$ is then naturally isomorphic to the functor
sending $$ M_{\coring{C}}\longrightarrow \prod_{x \in \,G} M^{cov(x)}, 
\text{ and }f \longrightarrow \left(\underset{}{}f_{M^{cov(x)}}\right)_{x \in
\,G},$$ where, for every $x \in  G$, $f_{M^{cov(x)}}$ is the
restriction of $f$ to $M^{cov(x)}$.
\end{remark}

Recall, from \cite{Nastasescu:1989a,Abrams/Menini:1999}, the coinduction functor
$\Coind{x}: \rmod{A_{\Sf{e}}} \rightarrow \grmod{A}$
which is given, over objects, by $$\Coind{x}(N)=\bigoplus_{y \in
\,G}\hom{A_{\Sf{e}}}{A_{y^{-1}x}}{N},$$ for every right
$A_{\Sf{e}}$-module $N$. The following shows  that the functor $\Coind{x}$ coincides, up to the isomorphism $\rcomod{\coring{C}} \cong \grmod{A}$, with the coinduction functor associted to the finitely generated projective right $\coring{C}$-comodule $[x]A$. Of course this involves the rational
functor $\Rat^{\T}$ constructed in Proposition \ref{smashProduct}.

\begin{proposition}\label{Coind}
Let $A$ be a $G$-graded ring by a group $G$ with neutral element $\Sf{e}$, $\coring{C}
= AG$ it associated $A$-coring, and $\T=(\coring{C},B,\esc{-}{-})$
the rational pairing of Proposition \ref{smashProduct}(ii).
Consider for some $x \in G$ the right $\coring{C}$-comodule $[x]A$ and its right
$B$-dual $([x]A)^{\star}=\hom{B}{[x]A}{B}$ as a
$(B,A_{\Sf{e}})$-bimodule. Then the functor $\Rat^{\T} \, \circ \,
\hom{A_{\Sf{e}}}{([x]A)^{\star}}{-}: \rmod{A_{\Sf{e}}}
\rightarrow \rcomod{\coring{C}}$ is right adjoint to the
coinvariant functor $\hom{\coring{C}}{[x]A}{-}:
\rcomod{\coring{C}} \rightarrow \rmod{A_{\Sf{e}}}$. Moreover, we have a
commutative diagram
$$\xymatrix@C=150pt@R=30pt{ \rmod{A_{\Sf{e}}}
\ar@{->}^-{\Rat^{\T} \, \circ \,
\hom{A_{\Sf{e}}}{([x]A)^{\star}}{-}}[r]
\ar@{->}_-{\Coind{x}}[rd] & \rcomod{\coring{C}} \ar@{->}^-{\cong}[d] \\
& \grmod{A}. }$$
\end{proposition}
\begin{proof}
The first statement is a direct consequence of Propositions
\ref{AdjointCoinv} and \ref{smashProduct}(iii), since we know by the isomorphism of equation \eqref{coinv-x} that $[x]A$
is finitely generated and projective right $\coring{C}$-comodule.

Let $N$ be any right $A_{\Sf{e}}$-module and $x \in G$. Then, we
have
\begin{eqnarray*}
  \Rat^{\T}(\hom{A_{\Sf{e}}}{([x]A)^{\star}}{N})  &=&
  \oplus_{y \in \,G}
\hom{A_{\Sf{e}}}{([x]A)^{\star}}{N}e_y, \quad \text{ by Proposition \ref{smashProduct}(iii)} \\
   & \cong & \oplus_{y \in \,G} \hom{B}{e_yB}{\hom{A_{\Sf{e}}}{([x]A)^{\star}}{N}} \\
   & \cong & \oplus_{y \in \,G} \hom{A_{\Sf{e}}}{e_yB\tensor{B}([x]A)^{\star}}{N}  \\
   & \cong & \oplus_{y \in \,G}
   \hom{A_{\Sf{e}}}{e_y([x]A)^{\star}}{N}.
\end{eqnarray*} On the other hand, one can prove that the mutually inverse maps
$$\xymatrix@R=0pt{  & e_y\left( [x]A \right)^{\star} \ar@{->}[rr] & & A_{y^{-1}x} &
\\ \left[ e_y\sigma \right. \ar@{|->}[r] & \left.\esc{y}{\sigma(1)}\right] &  &
\left[ \,\,[a \mapsto e_y\,\td{a_{y^{-1}x} a}] \right. & \ar@{|->}[l]
\left. a_{y^{-1}x}\,\, \right] }$$ establish an isomorphism of
right $A_{\Sf{e}}$-modules. Therefore, we obtain a natural isomorphism
$$\Rat^{\T}(\hom{A_{\Sf{e}}}{([x]A)^{\star}}{N}) \,\, \cong
\,\, \oplus_{y \in \,G}\hom{A_{\Sf{e}}}{A_{y^{-1}x}}{N}\,\, = \,\,
\Coind{x}(N).$$ This finishes the proof since the
compatibility, in relation with arrows, is clear.
\end{proof}

Let $x \in G$ and consider its associated right $\coring{C}$-comodule $[x]A$. Denote by $\scr{C}_x$ and $\scr{T}_x$, respectively, the torsion-free class and torsion class attached to this comodule, see Section \ref{Sec-2}. Up to the isomorphism of categories $\rcomod{\coring{C}} \cong \grmod{A}$, it is clear that a $G$-graded module $M$ belongs to $\scr{C}_x$ if and only if it $x$-homogeneous component vanishes. That is, $M \in \scr{C}_x$ if and only if $M_x=0$. A more specific computation of the torsion part of any $G$-graded right module can be given as follows. First denote by $\fk{r}_x$ the associted idempotent radical as in equation \eqref{radical}. 

\begin{lemma}\label{gr-radical}
Let $A$ be a graded ring by a group $G$ with neutral element $\Sf{e}$. Then, for every element $x \in G$ and $G$-graded right $A$-module $M$, we have 
$$\fk{r}_x(M)\,\,=\,\, {\rm Im}\lr{\oplus_{y \in\, G} \lr{M_x \tensor{A_{\Sf{e}}} A_{x^{-1}y}} \longrightarrow M},$$ wherein the map is the graded morphism given by the right action of $A$.
\end{lemma}
\begin{proof}
Fix an element $x$ in the group $G$. Let us first check that $\oplus_{y \in G} (M_x\tensor{A_{\Sf{e}}}A_{x^{-1}y})$ is an object in $\scr{T}_x$, for every $G$-graded right $A$-module $M$. So given a $G$-graded right module $Y$ in $\scr{C}_x$, and  graded map $f: \oplus_{y \in G} (M_x\tensor{A_{\Sf{e}}}A_{x^{-1}y}) \to Y$, we need to show that $f=0$. Since $Y_x=0$, it is obvious that $f(M_x \tensor{A_{\Sf{e}}}A_{\Sf{e}}) =0$. Taking an element in $M_x\tensor{A_{\Sf{e}}}A_{x^{-1}y}$ of the form $m \tensor{A_{\Sf{e}}}a$, we can find a finite sets $\{a_i\} \subset A_{\Sf{e}}$ and $\{a_i'\} \subset A_{x^{-1}y}$ such that $m  \tensor{A_{\Sf{e}}} a = \sum_i m\tensor{A_{\Sf{e}}}a_i a_i'$. Thus 
$$ f(m  \tensor{A_{\Sf{e}}} a)\, =\, \sum_i f( m\tensor{A_{\Sf{e}}}a_i) a_i' \,=\,0,$$ 
by linearity, which implies that $f=0$. Now, take an object $Z \in \scr{T}_x$, and a graded morphism $g : Z \to M$. We know that the canonical projection $Z \to Z/\oplus_{y \, \in G}(Z_x \tensor{A_{\Sf{e}}} A_{x^{-1}y})$ is zero since $Z/\oplus_{y \, \in G}(Z_x \tensor{A_{\Sf{e}}} A_{x^{-1}y}) \in \scr{C}_x$. That is, $Z \cong \oplus_{y \, \in G}(Z_x \tensor{A_{\Sf{e}}} A_{x^{-1}y})$. It clear then that $g$ factors thoughout the canonical injection ${\rm Im}(\nu) \hookrightarrow M$, where $\nu: \oplus_{y \, \in G}\lr{M_x \tensor{A_{\Sf{e}}} A_{x^{-1}y}} \to M $ is the obvious graded map, and this finishes the proof.
\end{proof}

Let us denote by $\scr{S}$ the set of isomorphisms classes $[Y]$ represented by simple right $A_{\Sf{e}}$-modules $Y$. $\scr{S}_{x}$ will denotes the set of isomorphisms classes $[S]$ represented by simple $G$-graded right $A$-module $S$ such that $ x \in \Sop{S}$.

The following  is a direct consequence of Corollary \ref{rep-simples}, Propositions \ref{smashProduct}, \ref{Coind}, and Lemma \ref{gr-radical} (compare with \cite{Menini/Nastasescu:1990}). 

\begin{corollary}\label{gr-simple}
Let $A$ be a graded ring by a group $G$. For every element $x \in G$, there is a bijective map 
$$
\xymatrix@R=0pt@C=40pt{ \scr{S}_{x} \ar@{->}[rr] & & \scr{S} \\ [S] \ar@{|->}[rr] & & [S_x] \\ \left[\underset{}{} \oplus_{y\in \, G}(Y\tensor{A_{\Sf{e}}}A_{x^{-1}y})\right]  & &  \ar@{|->}[ll] [Y]. }
$$ In particular every simple $G$-graded right $A$-module is semisimple right $A_{\Sf{e}}$-module.
\end{corollary}

\begin{remark}\label{gr-A.Obs}
The results presented in Section \ref{Sec-2} can be also applied to semi-group graded modules \cite{Abrams/Menini:1999}, $G$-set graded modules \cite{Nastasescu/Raianu/Oystaeyen:1990}, or more general to any category of entwined modules \cite{Brzezinski/Wisbauer:2003}. For instance, Corollary \ref{rep-simples} in conjunction with semigroup-graded version of Propositions \ref{smashProduct}, \ref{Coind}, will leads to \cite[Theorem 2.8]{Abrams/Menini:1999}.
\end{remark}

\providecommand{\bysame}{\leavevmode\hbox
to3em{\hrulefill}\thinspace}
\providecommand{\MR}{\relax\ifhmode\unskip\space\fi MR }
\providecommand{\MRhref}[2]{%
  \href{http://www.ams.org/mathscinet-getitem?mr=#1}{#2}
} \providecommand{\href}[2]{#2}

\end{document}